\documentclass[preprint,10pt]{elsarticle}
\usepackage[T1]{fontenc}
\usepackage{bm, amsmath, amsthm, amscd, amssymb, fancyhdr, upgreek, amssymb, subcaption, hyperref, enumerate, comment, siunitx, enumitem, graphicx, mathrsfs, mathtools}
\usepackage[dvipsnames]{xcolor}
\usepackage{microtype}
\usepackage[toc]{appendix}
\usepackage[parfill]{parskip}
\usepackage[all,2cell,ps]{xy}

\usepackage{hyperref}
\usepackage{url}
\usepackage{color}
\usepackage{breakurl}
\newcommand{\bburl}[1]{\textcolor{blue}{\url{#1}}}
\newcommand{\js}[1]{{#1\overwithdelims () p}}
\newcommand{\zsum}[1]{ \sum_{#1 = 0}^{p-1} }

\newcommand{\twocase}[5]{#1 \begin{cases} #2 & \text{#3}\\ #4
&\text{#5} \end{cases}   }

\newtheorem{theorem}{Theorem}[section]

\newtheorem{lemma}[theorem]{Lemma}

\theoremstyle{remark}
\newtheorem*{remark}{Remark}

\DeclareMathOperator{\rank}{rank}

\renewcommand\mod{~\mathrm{mod}~}

\usepackage{tikz}
\usepackage{tkz-tab}
\usepackage{tkz-graph}
\usetikzlibrary{shapes.geometric,positioning}

\numberwithin{equation}{section}

\newcommand{\hr}[1]{\href{#1}{\url{#1}}}

\journal{Conference Proceedings of ICERM}
\date{\today}

\begin{document}

\begin{frontmatter}
\title{Applications of Moments of Dirichlet Coefficients in Elliptic Curve Families
}

\affiliation[1]{organization={Pomona College}, addressline={Department of Mathematics and Statistics, 610 North College Avenue}, city={Claremont}, postcode={91711}, state={CA}, country={USA}}
\affiliation[2]{organization={University of Cambridge}, addressline={Department of Pure Mathematics and Mathematical Statistics}, city={Cambridge}, postcode={CB3}, state={0WA}, country={UK}}
\affiliation[3]{organization={University of Michigan}, addressline={Department of Mathematics}, city={Ann Arbor}, postcode={48104}, state={MI}, country={USA}}
\affiliation[4]{organization={Williams College}, addressline={Department of Mathematics and Statistics}, city={Williamstown}, postcode={01267}, state={MA}, country={USA}}
\affiliation[5]{organization={Yale University}, addressline={Department of Mathematics, PO Box 208283}, city={New Haven}, postcode={06520-8283}, state={CT}, country={USA}}

\author[1]{Zo\"{e} Batterman} \ead{zxba2020@mymail.pomona.edu}
\author[2]{Aditya Jambhale} \ead{aj644@cam.ac.uk}
\author[4]{Steven J. Miller} \ead{sjm1@williams.com}
\author[3]{Akash L. Narayanan} \ead{anaray@umich.edu}
\author[2]{Kishan Sharma} 
\ead{kds43@cam.ac.uk}
\author[2]{Andrew Yang} 
\ead{aky30@cam.ac.uk}
\author[5]{Chris Yao} 
\ead{chris.yao@yale.edu}

\begin{abstract} 
The moments of the coefficients of elliptic curve $L$-functions are related to numerous important arithmetic problems. Rosen and Silverman proved a conjecture of Nagao relating the first moment of one-parameter families satisfying Tate's conjecture to the rank of the corresponding elliptic surface over $\mathbb{Q}(T)$; one can also construct families of moderate rank by finding families with large first moments. Michel proved that if $j(T)$ is not constant, then the second moment of the family is of size $p^2 + O(p^{3/2})$; these two moments show that for suitably small support the behavior of zeros near the central point agree with that of eigenvalues from random matrix ensembles, with the higher moments impacting the rate of convergence. 

In his thesis, Miller noticed a negative bias in the second moment of every one-parameter family of elliptic curves over $\mathbb{Q}$ whose second moment had a (by him) calculable closed-form expression, specifically \emph{the first lower order term which does not average to zero is on average negative}. This \emph{Bias Conjecture} has now been confirmed for many families; however, these are highly non-generic families as they are specially chosen so that the resulting Legendre sums can be determined. For cohomological reasons, each subsequent term in the second moment expansion is smaller than the previous by a factor on the order of $\sqrt{p}$, and thus numerically, it is hard to see a term of size $p$ with a small negative average as it can be masked by a term of size $p^{3/2}$ which averages to zero.

Inspired by the recent successes by Yang-Hui He, Kyu-Hwan Lee, Thomas Oliver, Alexey Pozdnyakov and others in investigations of murmurations of elliptic curve coefficients with machine learning techniques, we pose a similar problem for trying to understand the Bias Conjecture. As a start to this program, we numerically investigate the bias conjecture and provide a visual representation of the bias for the second moment. We find a one-parameter family of elliptic curves whose bias is positive for half the primes. However, the numerics do not offer conclusive evidence that negative bias for the other half is enough to overwhelm the positive bias. Without an explicit expansion for the second moment, we are not able to extract potential negative bias of the order $p$ term.
\end{abstract}

\begin{keyword} 
Elliptic Curves \sep Second Moments \sep Bias Conjecture
\MSC[2020] 11G05 \sep 11G20 \sep 11G25
\end{keyword}

\end{frontmatter}
\section{Introduction}

We\footnote{This work was completed during the 2023 SMALL REU program at Williams College. It was supported in part by NSF Grants DMS1561945 and DMS1659037, the University of Michigan, the Churchill Foundation, and Williams College. We thank the organizers and participants of the workshop at ICERM on murmurations, especially Noam Elkies and Adam Logan, for many helpful conversations.} assume the reader is familiar with the basics of elliptic curves; see for example \cite{Sil94, Sil09}. Let $\mathcal{E} \to \mathbb{P}^1$ be a (non-split) elliptic surface over $\mathbb{Q}$, with Weierstrass equation
\begin{align}
\mathcal{E}:y^2\ =\ x^3+A(T)x+B(T),
\end{align}
with 
$A(T),B(T) \in \mathbb{Z}(T)$ and $4A(T)^3 +27B(T)^2 \neq 0$, with $\mathcal{E}$ a rational surface if $0 < \max \{3\deg A(T),2\deg B(T)\}<12$. Specializing $T$ to integers $t$, we see for all but finitely many choices we obtain an elliptic curve, which we denote by $E_t$:
\begin{align}
E_t \ \coloneqq \ y^2 = x^3 + A(t)x + B(t);
\end{align}
we thus have a one-parameter family of elliptic curves over  $\mathbb{Q}$. We consider the trace of Frobenius $a_t(p) := p + 1 - \#E_t(\mathbb{F}_p)$ where $\#E_t(\mathbb{F}_p)$ is the number of solutions of $E_t$ mod $p$ including the point at infinity.
By the expression of $E_t$, we have an explicit formula for $a_t(p)$ as a sum of Legendre symbols $\js{a}$, which is 1 if $a$ is a non-zero square modulo the prime $p$, 0 if $a$ is zero modulo $p$, and -1 otherwise:
\begin{align}
- a_t(p) \ = \ \sum_{x = 0}^{p-1} \left( \frac{x^3 + A(t)x + B(t) }{p} \right).
\end{align}
There are simple closed form expressions for the sum of $\js{f(t)}$ over $t$ modulo $p$ if $f$ is linear or quadratic; see for example Appendix A of \cite{AL-RM07}.

\begin{lemma}[Linear and Quadratic Legendre Sums]\label{labquadlegsum} Assume $p > 2$. If $a$ is not zero modulo $p$ then \begin{equation} \zsum{t} \js{at + b} \ = \ 0, \end{equation} and if 
$a$ and $b$ are not both zero modulo $p$ then
\begin{equation}
\twocase{\zsum{t} \js{at^2 + bt + c} \ = \ }{(p-1)\js{a}}{{\rm if $p
| (b^2 - 4ac)$}}{-\js{a}}{{\rm otherwise}.}
\end{equation} 
\end{lemma}

Define the $m$\textsuperscript{th} \emph{moment} (we do not normalize by $1/p$) of $\mathcal{E}$ to be
\begin{align}
\mathcal{A}_{m,\mathcal{E}}(p) \ := \ \sum_{t = 0}^{p-1} a_t(p)^m.
\end{align} 
These moments play a key role in determining the arithmetic properties of elliptic curves. Rosen and Silverman \cite{RS98} proved a conjecture of Nagao \cite{Nag97}, which states that if Tate's conjecture holds, then
\begin{align}
\lim_{X\to \infty} \frac{1}{X} \sum_{p\leq X} \frac{\mathcal{A}_{1,\mathcal{E}}(p)\log p}{p} \ = \ -\rank \mathcal{E}(\mathbb{Q}(T));
\end{align} thus the rank of the elliptic surface is determined by the first moment when Tate's conjecture is true. In particular, Tate's conjecture is known to hold for rational surfaces \cite{Shi72}. Equivalently, the above limit tells us there is a negative bias in the coefficients $a_t(p)$, and that bias is related to the rank with larger biases yielding greater ranks. 

Turning to the second moment, Michel proved the following.

\begin{theorem}[Michel \cite{Mic95}]
For an elliptic surface $\mathcal{E}$ with non-constant $j(T)$-invariant, the second moment is of the form
\begin{align}
\mathcal{A}_{2,\mathcal{E}}(p) \ = \ p^2 + O(p^{3/2}).\end{align} In particular, the lower order terms arise from cohomological arguments, and each is a factor on the order of $\sqrt{p}$ less than the previous.
\end{theorem}

These two moments suffice to show agreement between zeros near the central point, i.e. the low-lying zeros, and eigenvalues of the corresponding random matrix ensemble \cite{Mil04}, providing support for the Katz-Sarnak density conjectures. The higher moments do not affect the main terms; similar to the Berry-Essen theorem on convergence in the Central Limit Theorem, these only impact the rate of convergence.

Recall the negative bias found for the first moment. The natural question to consider is whether higher moments share the same negative bias. In his thesis, Miller \cite{Mil02} noted that for every family where he could find a closed form expression for the second moment, the first lower order term which did not average to zero had a negative average.\footnote{These families were such that upon switching the order of summation one obtained a Legendre sum of a polynomial that was linear or quadratic in the variable of interest, which meant the resulting sum had a simple closed form expression; by carefully choosing the family the remaining sums were also computable.} For a generic family, we do not expect to find such a closed form expression using the Legendre sum approach. This is due to the general intractability of cubic and higher Legendre sums, which explains the difficulty of working with elliptic curve coefficients. 

Similar to how the negative bias in the first moment has important arithmetic applications to the average rank, a negative bias for the second moment has applications in understanding low-lying zeros: it results in a lower order correction term in the 1-level density which, for finite conductors, increases the bound for the average rank (and perhaps can help explain the elevated ranks one observes in families with small conductors); see \cite{Mil05} for details, which we briefly summarize for completeness in Appendix \ref{sct:applications}. Since then, the Bias Conjecture has been confirmed in every family studied \cite{A+23,KN21,KN22,Mil02,Mil04,Mil05,MMRW16}; however, these families are highly non-generic as they all have second moments with closed form expressions. 

Explicitly, the \emph{Bias Conjecture}, in its original formulation, states that the largest term in the second moment expansion which does not average to zero is on average negative. This is a weak formulation, and allows there to be a positive behavior for some, or even infinitely many, primes, so long as the other primes are more negative than those that are positive. We shall call this the \emph{Weak Bias Conjecture}, and by the \emph{Strong Bias Conjecture}, we mean that the largest lower order term in the second moment expansion which does not average to zero is negative except for finitely many $p$.

Inspired by recent remarkable successes by Cowan, Yang-Hui He, Kyu-Hwan Lee, Thomas Oliver, Alexey Pozdnyakov and others \cite{Co23, HLOa, HLOb, HLOc, HLOP22} in investigating murmurations of elliptic curve coefficients with machine learning techniques, we pose a similar problem to understand the Bias Conjecture. Thus, we search families of elliptic curves where we can no longer compute in closed form the second moment, and we try to see if there is a negative bias on average in the first lower order term not averaging to zero. Computationally, this is a very difficult problem as a term of size $p$ that is on average negative is dwarfed by a term of size $p^{3/2}$ which averages to zero. We thus believe that the machinery and methods of recent murmuration papers could be successfully applied here.

As a start to such a program, below we report on numerical and theoretical investigations for a one-parameter family of elliptic curves. We show that for half of the primes there is a positive bias, thus disproving the strong form of the conjecture. The numerics are not convincing enough to determine if the bias conjecture holds. We hope that the machine learning techniques used for murmurations can tell us whether or not a negative bias from any remaining primes overwhelms the positive bias. If it does, then this would provide non-trivial support for the weak form of the conjecture.

\section{Known Biases.}
The table below lists first and second moments of several elliptic curve surfaces. We gather previous expressions for the second moment derived in \cite{A+23,Mil02,Mil05} in the following table.
Let $n_{3,2,p}$ equal to the number of cube roots of 2 modulo $p$, and set $c_0(p) = \left[ \left( \frac{-3}{p} \right) + \left( \frac{3}{p} \right)  \right] p$, $c_1(p) = \left[ \sum_{x = 0}^{p-1} \left( \frac{x^3-x}{p} \right)  \right] ^2$, $c_{3 /2}(p) = p \sum_{x = 0}^{p-1}\left( \frac{4x^3+1}{p} \right) $, and $c_{p,a;m}=1$ if $p$ congruent to $a\mod m$ and otherwise equal to $0$.

\begin{center}
\begin{tabular}{c c c}
\text{Family $\mathcal{E}$} & $A_{1,\mathcal{E}}(p)$ & $A_{2,\mathcal{E}}(p)$ \\
\hline 
$y^2 = x^3 + Sx + T$ & 0 & $p^3 + p^2$ \\
$y^2 = x^3+2^4(-3)^3(9T+1)^2$ & 0 & $\begin{cases} 2p^2 - 2p & \qquad p \equiv 2 \mod 3 \\
0 & \qquad p \equiv 1 \mod 3 \end{cases}$ \\
$y^2 = x^3\pm 4(4T+2)x$ & 0 & $\begin{cases} 2p^2 - 2p & \qquad p \equiv 1 \mod 4 \\
0 & \qquad p \equiv 3 \mod 4 \end{cases}$ \\
$y^2 = x^3 + (T+1)x^2 + Tx$ & 0 & $p^2 -2p-1$ \\
$y^2 = x^3+x^2+2T+1$ & 0 & $p^2 - 2p - \left(\frac{-3}{p}\right)$ \\
$y^2=x^3 + Tx^2 +1$ & $-p$ & $p^2 - n_{3,2,p}p-1+c_{3/2}(p)$ \\
$y^2=x^3-T^2x+T^2$&$-2p$&$p^2-p-c_1(p)$\\
$y^2=x^3-T^2x+T^4$&$-2p$&$p^2-p-c_1(p) -c_0(p)$\\
$y^2=x^3+Tx^2-(T+3)x+1$&$-2c_{p,1;4}p$&$p^2-4c_{p,1;6}p-1$
\end{tabular}
\end{center}

Note $c_1(p)$ is the square of the coefficients from an elliptic curve with complex multiplication. It is non-negative and of size $p$ for $p \not\equiv 3 \mod 4$ and size 0 for $p \equiv 1 \mod 4$ (send $x \mapsto -x \mod p$).
Except for the family $y^2 = x^3 + Tx^2 + 1$, all the surfaces $\mathcal{E}$ have $A_{2,\mathcal{E}}(p) = p^2 - h(p)p + O(1)$, where $h(p)$ is non-negative. It is somewhat remarkable that all of these families have a correction to the main term in Michel's theorem in the same direction, and we review the consequence this has on the average rank in \ref{sct:applications}; this application was the motivation behind the initial interest in the conjecture. 

The majority of surfaces in the table lie in the more general collection of cubic pencils $\mathcal{C}: y^2 = P(x)T + Q(x)$ where $P(x),Q(x)\in \mathbb{Z}[x]$ are degree at most 3. Kazalicki and Naskrecki \cite{KN21} prove the second moment of an elliptic surface in $\mathcal{C}$ satisfies the Bias Conjecture by counting points on Kummer threefolds. Since they resolved the more general case, they were thus able to recover the second moment for the linear elliptic curve families in \cite{A+23}: $y^2=(ax^2+bx+c)(dx+e+T)$, $y^2=x(ax^2+bx+c+dTx)$, and $y^2=x(ax+b)(cx+d+Tx)$.

We state results for those surfaces with higher degree than 2 in $T$.
In the sequel paper \cite{KN22}, they prove the Bias Conjecture holds for the second moment of the surface $y^2 = x(T^2(1+T^2)^3+2(1+T^2)^2x+x^2)$. Asada \emph{et al.} in \cite{A+23} explicitly found all moments for families with parameter $T$ of degree greater than 1: $y^2=x^3-(cT+d)^rAx$, $y^2=x^3+(cT+d)^rB$, and $y^2=x^3+T^2Ax+T^3B$ where $A,B,c,d \in \mathbb{Z}$.
We remark, however, that those families whose closed-form expression were explicitly found have coefficient polynomials with low degree. Since these families are so special, we may be seeing non-generic behavior, which motivates the search for the second moment of a family with larger degree polynomials in its definition.

\section{Numerical Positive Bias.}
We study the rational elliptic surface given by
\begin{equation}
\mathcal{F}: y^2 \ =\ x^3 + x + T^3.
\end{equation}
We show for this family that there is a positive bias for half the primes.


\begin{lemma}\label{lemma:second-moment-2-mod-3}
For primes congruent to 2 modulo 3, the second moment of $\mathcal{F}$ is given by
\begin{equation}
\mathcal{A}_{2,\mathcal{F}}(p) \ =  \ p^2 + p.
\end{equation}
\end{lemma}

\begin{proof}
Fix a prime $p$ congruent to 2 modulo 3. By definition, we write the second moment
\begin{align}
\mathcal{A}_{2,\mathcal{F}}(p)\ =\ \sum_{x = 0}^{p-1} \sum_{w = 0}^{p-1} \sum_{t = 0}^{p-1}\left(\frac{x^3+x+t^3}{p}\right)\left(\frac{w^3+w+t^3}{p}\right).
\end{align}
By the multiplicativity of the Legendre symbol, we get
{\small 
\begin{equation*}
\mathcal{A}_{2,\mathcal{F}}(p)\ =\ \sum_{x = 0}^{p-1}\sum_{w = 0}^{p-1}\sum_{t = 0}^{p-1}\left(\frac{t^6 + t^3(x^3 + x + w^3 + w) + x^3 w^3 + x^3 w + x w^3 + xw}{p}\right).
\end{equation*}
}
Since $p$ is congruent to 2 modulo 3, the map $t^3 \mapsto t$ is an automorphism which allows us to reduce the sextic in $t$ to a quadratic in $t$:
{\footnotesize
\begin{equation}\label{quadratic}
\mathcal{A}_{2,\mathcal{F}}(p) \ =\ \sum_{x = 0}^{p-1}\sum_{w = 0}^{p-1}\sum_{t = 0}^{p-1} \left(\frac{t^2 + t(x^3 + x + w^3 + w) + x^3 w^3 + x^3 w + x w^3 + xw.}{p}\right).
\end{equation}
}
The expression for the Legendre symbol of a quadratic is tractable, as remarked in Lemma \ref{labquadlegsum}. We consider the discriminant of the quadratic polynomial,
\begin{equation}
    \Delta(x,w) \ \coloneqq\ (x^3 + x + w^3 + w)^2 - 4(x^3 w^3 + x^3 w + x w^3 + xw),
\end{equation}
which simplifies to $\Delta(x,w) = (x^3 + x - w^3 - w)^2$.
We note $\Delta(x,w) = 0 \mod p$ if and only if $\delta(x,w) \coloneqq x^3 + x - w^3 - w \equiv 0 \mod p$,
if and only if
\begin{equation}
    x = w \quad \text{or} \quad \delta_2(x,w) \coloneqq x^2 + wx + w^2 + 1 = 0 \mod p.
\end{equation}
We wish to determine the number of solutions in each case. For $\delta_2 \equiv 0 \mod p$, we get from the quadratic formula modulo $p$ that
\begin{equation}
    x\ =\ \frac{-w\pm\sqrt{-3w^2-4}}{2}.
\end{equation}
There are
\begin{equation}
p + \left(\frac{-1}{p}\right) \sum_{w= 0}^{p-1} \left(\frac{3w^2+4}{p}\right)\  =\ p - \left(\frac{-1}{p}\right) \left(\frac{3}{p}\right)
\end{equation}
solutions to $\delta_2(x,w) = 0 \mod p$. 
If $p \equiv 11 \mod 12$, we have that $(3/p) = 1$, while $(-1/p) = -1$, leading to $(3/p)(-1/p)= -1$. On the other hand, if $p \equiv 5 \mod 12$, we have that $(3/p) = -1$, while $(-1/p) = 1$. Thus, we also get $(3/p)(-1/p)= -1$. In either case, there are $p+1$
solutions to $\delta_2(x,w) \equiv 0 \mod p$.

We also note that $x = w$ and $\delta_2(x,w) \equiv 0 \mod p$ cannot both hold. If both were true, then we would have $3x^2 + 1 \equiv 0 \mod p$,
implying that $-3$ is a square modulo $p$. However, by above work, we have $(-3/p) = -1$, a contradiction. The $\#\mathbb{F}_p = p$ solutions to $x = w$ are disjoint from the $p+1$ solutions $\delta_2(x,w) \equiv 0 \mod p$. Hence, there are $2p + 1$
solutions to $\Delta(x,w) \equiv 0 \mod p$.

Now, we may write \eqref{quadratic} as
\begin{align}
& \sum_{\substack{(x,w) \\ p \mid \Delta (x,y)}} \sum_{t = 0}^{p-1} \left(\frac{t^2 + t(x^3 + x + w^3 + w) + x^3 w^3 + x^3 w + x w^3 + xw}{p}\right) \\
&+ \sum_{\substack{(x,w) \\ p \nmid \Delta (x,y)}} \sum_{t = 0}^{p-1}  \left(\frac{t^2 + t(x^3 + x + w^3 + w) + x^3 w^3 + x^3 w + x w^3 + xw}{p}\right).
\end{align}
Applying the formula for quadratic Legendre sums, we get
\begin{align}
    \mathcal{A}_{2,\mathcal{F}}(p) &= (p-1) \sum_{\substack{(x,w) \\ \Delta(x,w) \equiv 0 \mod p }} \bigg(\frac{1}{p} \bigg) - \sum_{\substack{(x,w) \\ \Delta(x,w) \not\equiv 0 \mod p}} \bigg(\frac{1}{p} \bigg) \\&= (p-1) \cdot \#\{\Delta(x,w) \equiv 0 \mod p\} - \#\{\Delta(x,w) \not\equiv 0 \mod p\}.
\end{align}
As there are $2p+1$ solutions to $\Delta(x,w) \equiv 0 \mod p$ and thus $p^2 - 2p - 1$ solutions to $\Delta(x,w) \not\equiv 0 \mod p$, we get $\mathcal{A}_{2,\mathcal{F}}(p) = p^2 + p$, as desired.
\end{proof}
When $p$ is congruent to 1 modulo 3, the map $t^3\mapsto t$ fails to be an automorphism, and Lemma \ref{labquadlegsum} is no longer applicable. In this case, calculating the second moment explicitly cannot be done with standard techniques in number theory, and we instead turn to empirical evidence to gain a preliminary understanding of the second moment. Because our family has non-constant $j(T)$-invariant $j(T) = 4/(4+27T^3)$, we know 
\begin{equation}\label{eqn: michel second moment}
\mathcal{A}_{2, \mathcal{F}}(p) \ = \  p^2 + \gamma(p) p^{3/2} + \delta(p)p + O(\sqrt p),
\end{equation}
where $\gamma(p)$ and $\delta(p)$ are $O(1)$. We define the \emph{bias} of the second moment to be 
\begin{equation}\label{eqn: original bias def}
    B_{\mathcal{F}}(p) \ \coloneqq\  \mathcal{A}_{2,\mathcal{F}}(p) p^{-3/2} - p^{1/2}.
\end{equation} The rationale behind this definition is to subtract the known main term of $p^2$ and divide by the size of the largest possible lower order term for normalization purposes. Using \eqref{eqn: michel second moment}, we may rewrite the bias as 
\begin{equation}\label{biasDef}
B_{\mathcal{F}}(p)\ =\ \gamma(p) + \delta(p)p^{-1/2} + O(p^{-1}).
\end{equation} 

\begin{figure}[h]
    \centering
\includegraphics[scale=0.35]{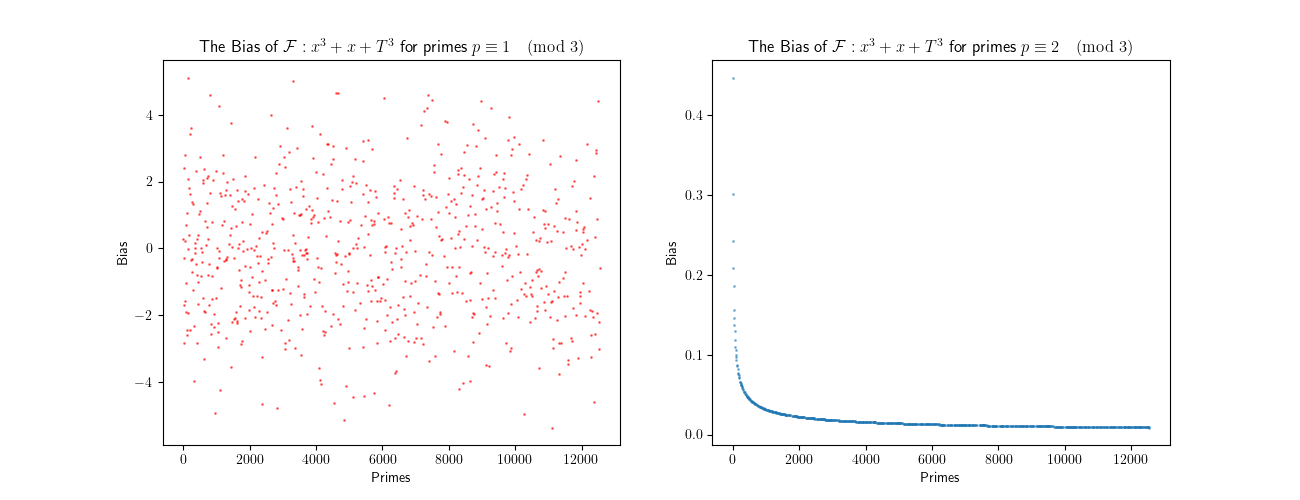}
    \caption{Left: A plot of the bias in the second moment for primes congruent to 1 mod 3. Right: The same plot but for primes congruent to 2 mod 3.}
    \label{fig:bias split by congruence}
\end{figure}

\begin{figure}[h]
    \centering
\includegraphics[scale=0.35]{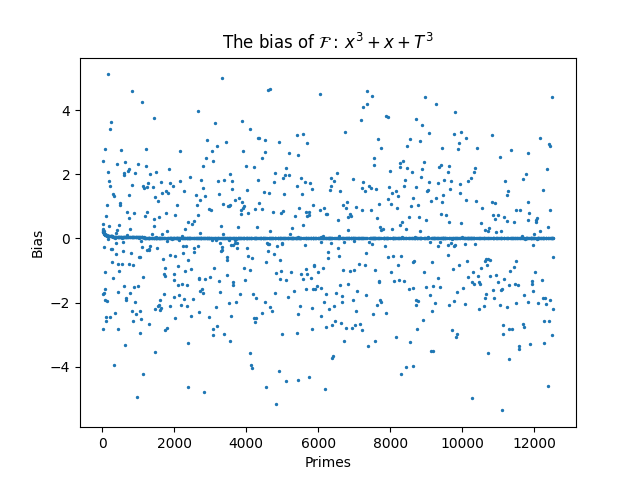}
\includegraphics[scale=0.35]{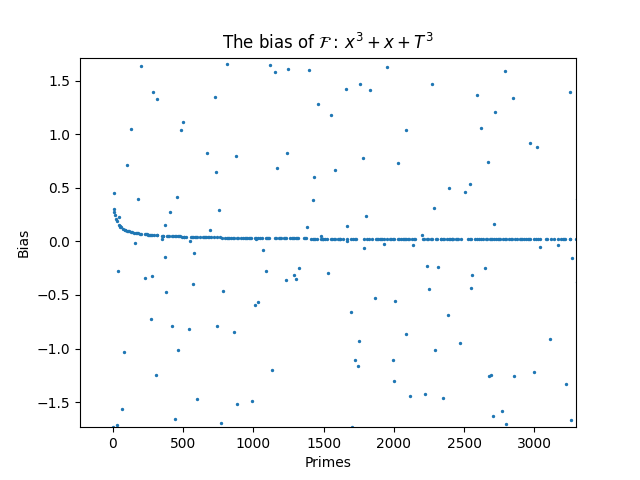}
    \caption{Left: A plot of the bias in the second moment for primes up to around 16500. Right: A zoomed-in version of the previous plot.}
    \label{fig:pos-bias}
\end{figure}

We obtained the second moment $\mathcal{A}_{2,\mathcal{F}}(p)$ numerically using code in C++ which may be accessed in \cite{Yao24}, and Figure \ref{fig:pos-bias} provides plots of the bias in the second moment for primes $p$ up to 16500. We note the presence of a pronounced ``tail" starting from the positive $y$-axis that appears to asymptotically approach zero. This phenomenon can be explained if we split the bias for primes in a residue class modulo 3, as in Figure \ref{fig:bias split by congruence}. By Lemma \ref{lemma:second-moment-2-mod-3}, for primes $p$ congruent to 2 modulo 3 we have $\delta(p) = 1$ and $\gamma(p) = 0$, and plugging this into \eqref{biasDef} gives that the bias is $B_\mathcal{F}(p) = p^{-1/2}$, which matches the tail in Figures \ref{fig:bias split by congruence} and \ref{fig:pos-bias}.

\begin{figure}[h]
    \centering
    \includegraphics[scale=0.38]{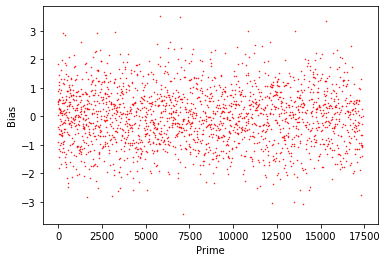}
    \includegraphics[scale=0.38]{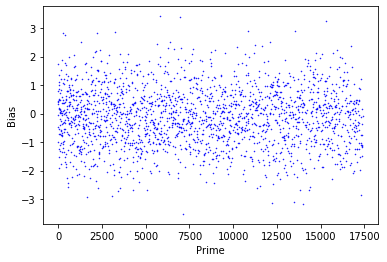}
    \caption{Plots of the bias for unbiased random second moments (red) versus biased random second moments (blue) for a random simulation.}
    \label{fig:simulated2ndMoment}
\end{figure}

In contrast, there appears to be no discernible pattern in the bias for primes congruent to 1 modulo 3. Motivated by the observation that the bias seems random, we simulated what the bias would look like if this were the case: for a given prime $p$, we generated $p$ many random moments according to the Sato-Tate distribution, i.e., $a_{t}(p) = 2\sqrt{p}\cos(\Theta)$
where $\Theta$ has probability density function $2\sin^2(\theta)/\pi$ for $\theta \in [0,\pi]$. Next, we summed across the squares of these random moments, simulating a random second moment for the one-parameter family $\mathcal{F}$; that is, a random simulation of $\mathcal{A}_{2,\mathcal{F}}(p)$.  As a comparison, we also simulated the second moment again, this time creating an artificial negative bias by subtracting a small $p^{3/2}$ term, i.e., $\mathcal{A}'_{2,\mathcal F}(p) = \mathcal{A}_{2,\mathcal F}(p) - 0.1\cdot p^{3/2}$, where $\mathcal{A}_{2,\mathcal F}(p)$ is the simulated (unbiased) second moment. In both cases, we calculated the bias according to \eqref{eqn: original bias def}, the results of which are shown in Figure \ref{fig:simulated2ndMoment}. The two graphs appear indistinguishable, meaning we cannot reasonably understand the average $p^{3/2}$ coefficient of our family simply from the graph of the bias. Instead, we turn to the \textit{average bias}:
\begin{equation}
    \lim_{x\to\infty} \frac{1}{\pi(x)}\sum_{p\leq x} B_{\mathcal{F}}(p) \ = \ \lim_{x\to\infty} \frac{1}{\pi(x)}\sum_{p\leq x}\big( \gamma(p) + \delta(p)p^{-1/2} + O(p^{-1})\big).
\end{equation}

\begin{figure}[h]
    \centering
    \includegraphics[scale=0.35]{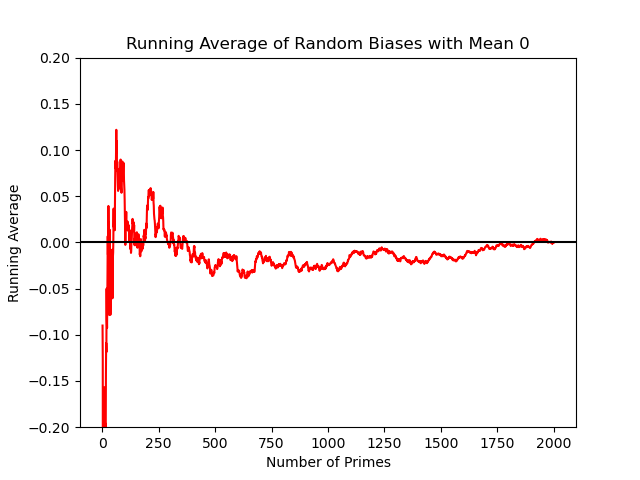}
    \includegraphics[scale=0.35]{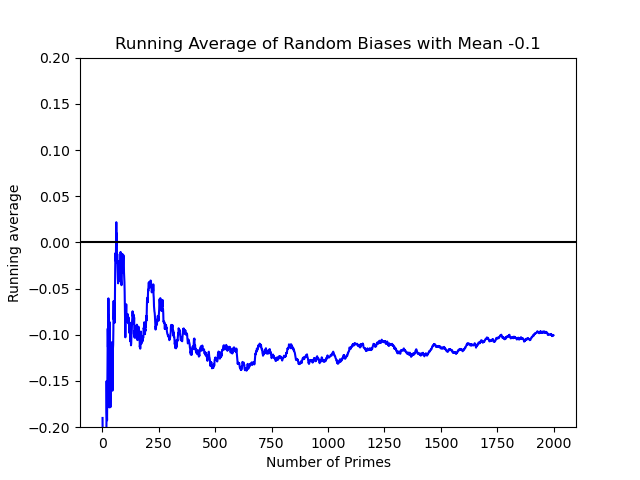}
    \caption{Unbiased running averages (red) versus biased running averages (blue) for a random simulation (Zoomed in)}
\label{fig:simulatedzoomedInRunningAverage}
\end{figure}

\begin{figure}[h]
\centering
\includegraphics[scale=0.35]{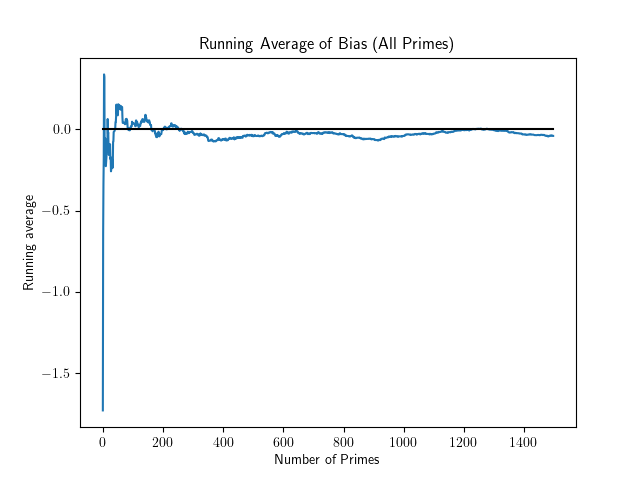}
\includegraphics[scale=0.35]{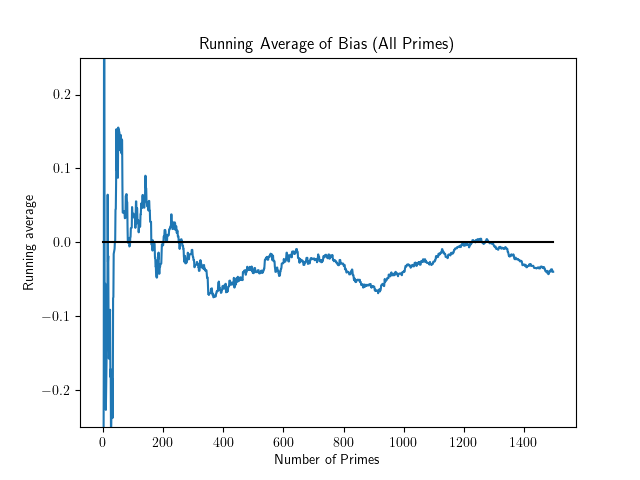}
\caption{Left: Running average of the bias for $\mathcal{F}: y^2 = x^3 + x + T^3$. Right: A zoomed-in version of the previous plot.}\label{fig:positiveBiasFamilyRunningAverage}
\end{figure}

For $x$ large, we have $\sum_{p\leq x} 1/p \sim \log\log(x)$ so $\pi(x)$ asymptotically dominates $\sum_{p\leq x} 1/p$, and we can drop the term from the limit, getting
\begin{equation}\label{biasaverage}
\lim_{x\to\infty} \frac{1}{\pi(x)}\sum_{p\leq x}B_{\mathcal{F}}(p) \ = \ \lim_{x\to\infty} \frac{1}{\pi(x)}\sum_{p\leq x}\big( \gamma(p) + \delta(p)p^{-1/2})\big)
\end{equation}
Next, using the integral bound on a sequence, we get
\begin{equation}
    \sum_{p\leq x} p^{-1/2} \ \leq \ \int_{1}^x t^{-1/2} dt \ = \ 2\sqrt{x} - 2 \ =  \ O(\sqrt{x}).
\end{equation}

\begin{figure}[h]
\centering
\includegraphics[scale=0.35]{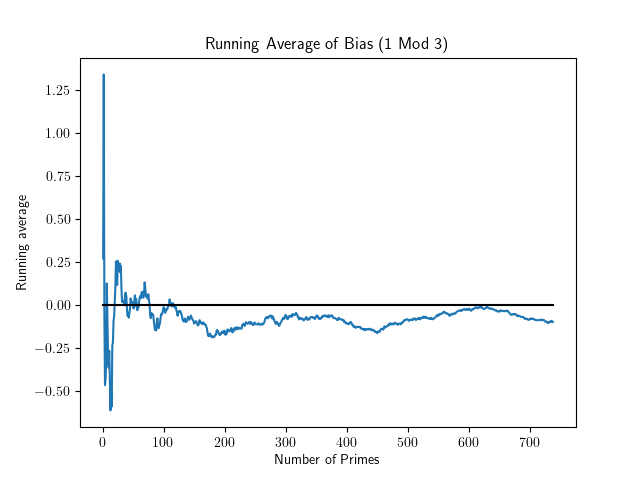}
\includegraphics[scale=0.35]{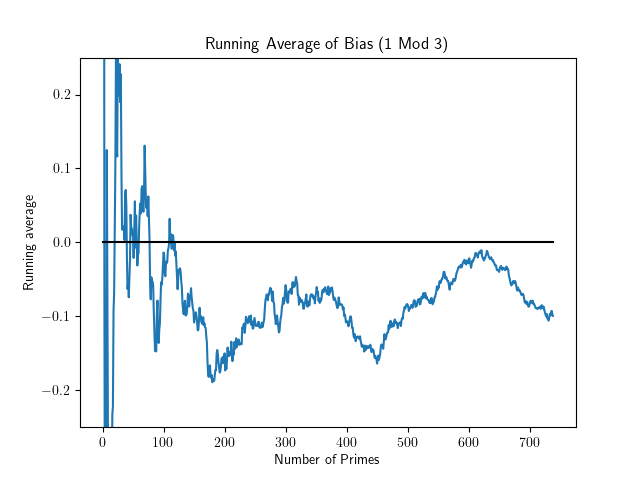}
\caption{Left: Running average of the bias for $\mathcal{F}: y^2 = x^3 + x + T^3$ for $p\equiv 1\mod 3 $. Right: A zoomed-in version of the previous plot.}\label{fig:positiveBiasFamilyRunningAverage 1 mod 3}
\end{figure}

Since $\delta(p)$ is $O(1)$, we have $\sum_{p\leq x} \delta(p) p^{-1/2} = O(\sqrt{x})$.
When we normalize this error term by $\pi(x)\sim x/\log(x)$, we see that it once again vanishes in the limit. Plugging the error term back into \eqref{biasaverage}, we obtain
\begin{equation}
    \lim_{x\to\infty} \frac{1}{\pi(x)}\sum_{p\leq x} B_{\mathcal{F}}(p)\ = \
    \lim_{x\to\infty} \frac{1}{\pi(x)}\sum_{p\leq x} \gamma(p).
\end{equation}
\begin{remark}
When we assume Michel's theorem has error term $O(p^a)$ for some $a < 3/2$, we obtain
\begin{equation}
\frac{1}{\pi(x)}\sum_{p\leq x} B_{\mathcal{F}}(p) \ = \ \frac{1}{\pi(x)}\sum_{p\leq x} \gamma(p) + O(x^{a-3/2}\log(x)).
\end{equation}
\end{remark}

\begin{figure}[h]
\centering
\includegraphics[scale=0.35]{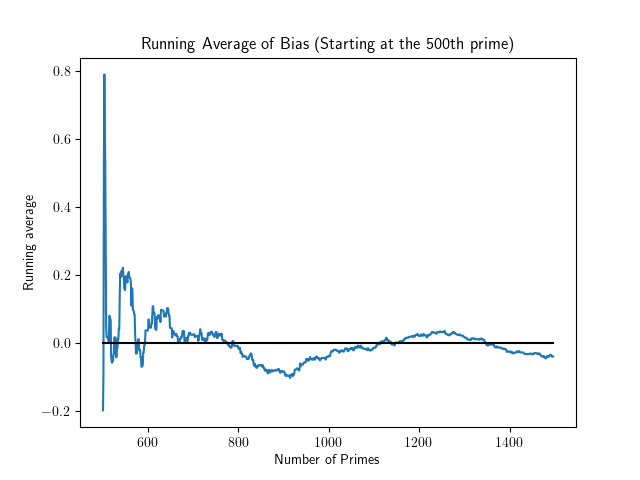}
\includegraphics[scale=0.35]{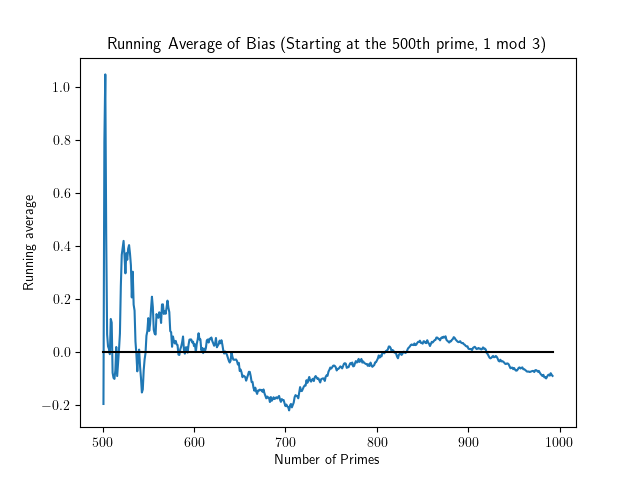}
\caption{Left: Running average of the bias for $\mathcal{F}: y^2 = x^3 + x + T^3$ for all primes, starting at the 500th prime. Right: Running average of the bias for $\mathcal{F}: y^2 = x^3 + x + T^3$ for primes $p \equiv 1\mod3$ starting at the 500th prime.}\label{fig:excised_data_our_family}
\end{figure}

Theoretically, of course, we may require large primes to observe the convergence if the constant in the big-O term is very large. As such, we again turned to our simulation of the second moment, both with and without the artificial bias, and calculated the average bias. Figure \ref{fig:simulatedzoomedInRunningAverage} show the running average bias (where the $x$-axis represents the number of primes summed over) for both the unbiased and biased random second moments. 
While the bias is hardly detectable in Figure \ref{fig:simulated2ndMoment}, it is quite clear in Figure \ref{fig:simulatedzoomedInRunningAverage}. This gives us the heuristic that the running average being close to $0$ at approximately $2000$ primes averaged over suggests that $\gamma(p)$ averages to 0. 

\begin{figure}
\centering
\includegraphics[scale=0.35]{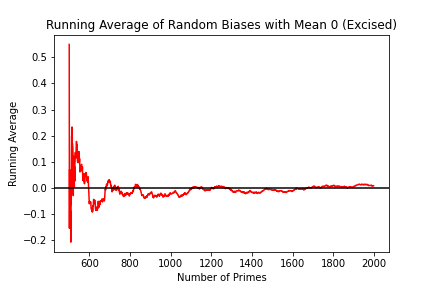}
\includegraphics[scale=0.35]{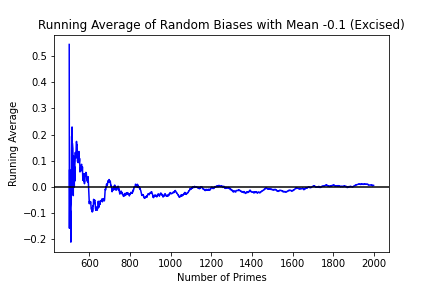}
\caption{Left: Unbiased running average (red) starting at the 500th prime versus biased running average (blue) starting at the 500th prime.}\label{fig:simulatedExcisedRunningAvg}
\end{figure}

Returning to our family, Figures \ref{fig:positiveBiasFamilyRunningAverage} and \ref{fig:positiveBiasFamilyRunningAverage 1 mod 3} show the running average of the bias for all primes and for primes congruent to 1 modulo 3, respectively. To prevent the effect of small primes, we also calculate the running averages starting at the 500th prime (e.g., summing over 100 primes would be the 500-600th primes). Figure \ref{fig:excised_data_our_family} shows this excised running average for our family, while Figure \ref{fig:simulatedExcisedRunningAvg} shows the excised running average for our simulation. Since the running average of the bias for the family in Figures \ref{fig:positiveBiasFamilyRunningAverage}, \ref{fig:positiveBiasFamilyRunningAverage 1 mod 3}, and \ref{fig:excised_data_our_family} mimics the unbiased running average of Figures \ref{fig:simulatedzoomedInRunningAverage} and \ref{fig:simulatedExcisedRunningAvg},  $\gamma(p)$ for our family appears to average to 0 for primes congruent to 1 modulo 3 and in general for all primes.

\begin{figure}[h]
\centering
\includegraphics[scale=0.35]{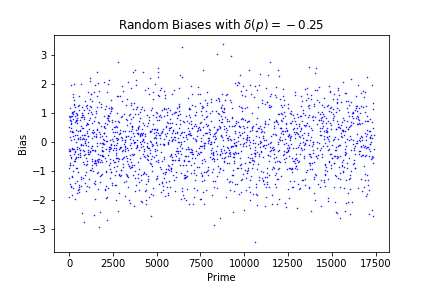}
\includegraphics[scale=0.35]{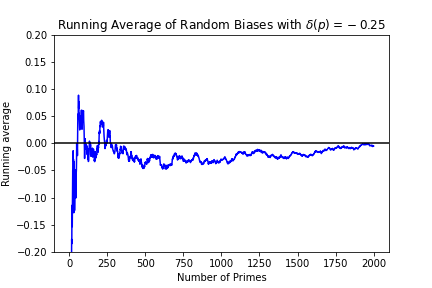}
\includegraphics[scale=0.35]{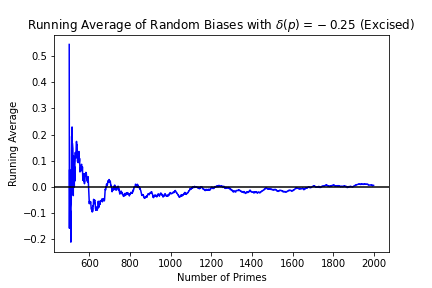}
\caption{Left: The bias for a random simulation with negative lower-order bias. Right: Running average of negative lower-order bias. Bottom: Running average of negative lower-order bias starting at the 500th prime.}\label{fig:lowerOrderBias}
\end{figure}

However, these numerics would not be able to distinguish if the family has lower-order bias, i.e., if $\mathcal{A}_{2, \mathcal F}(p) = p^2 + \gamma(p) p^{3/2} + cp$ where $\gamma(p)$ averages to 0 but $c\neq 0$. In this case, our bias would be $B_\mathcal{F}(p) = \gamma(p) + c/\sqrt p$, which would still average to 0, but the noise from the higher order term $\gamma(p)$ would also prevent us from isolating this lower order contribution. Figure \ref{fig:lowerOrderBias} depicts the results of a simulation where $\gamma(p)$ averages to 0 as before, but a lower order bias $\delta(p) = -0.25$ is introduced. As a result, we cannot conclude whether there is negative bias in the $p\equiv 1\mod{3}$ case or whether any negative bias in this case overwhelms the positive bias when $p\equiv 1\mod{3}$; hence, the Weak Bias Conjecture remains unresolved for our family.



\appendix
\section{Application to Excess Rank}\label{sct:applications}

We review from \cite{Mil02, Mil05} the connection between a negative bias in the second moment and the distribution of low-lying zeros.

It is possible to predict the average rank for a one-parameter family of elliptic curves. By assuming the Birch and Swinnerton-Dyer conjecture, Silverman's specialization theorem implies that eventually all curves $E_t$ have rank at least $r$ if our elliptic curve over $\mathbb{Q}(T)$ has rank $r$. In many families we expect half the curves to have even rank and half to have odd. In that case, random matrix theory predicts that, in the limit of large conductor, 50\% of elliptic curves in the family should be rank $r$ and 50\% should have rank $r+1$ for an average rank of $r + 1/2$.
However, a disagreement was noticed in many families between the predicted average rank and the observed average rank, which is called the \emph{excess rank problem}. The observed excess rank could easily be a result of a slow rate of convergence, as it is likely that the governing quantity is the logarithm of the conductor. Eventually, Watkins \cite{Wa07} went far enough to see the observed excess rank noticeably decreases. 

Miller noticed that a lower order negative bias increases the bound for average ranks in families through statistics of zero densities near the central point (see \cite{Mil05} for details). Unfortunately, this change happens by only a small amount which is not enough to explain the observed excess rank, but does motivate the study of these biases. 


We summarize the calculations from \cite{Mil05} to bound, under the GRH and the Density Conjecture, the average rank of a one-parameter family $\mathcal{E}$ of rank $r$ over $\mathbb{Q}(T)$. Letting $N$ be a positive integer, we study elliptic curves $E_t$ for $t \in [N, 2N]$; we then send $N \to \infty$. Let $C_t$ denote the conductor of an elliptic curve $E_t$, and thus the average spacing between zeros near the central point is on the order of $1/\log C_t$. Set $\log R = \frac{1}{N} \sum_{t=N}^{2N} \log C_t$, the average log-conductor (with more work one an normalize each curve's zeros by the correct local quantity).

To consider the effect of an average negative bias in the second moment, we assume that the second moment is of the form $p^2 - m_{\mathcal{{E}}}p + O(1)$ where $m_{\mathcal{E}} > 0$. For $\phi$ an even Schwartz test function supported on the interval $(-\sigma, \sigma)$, one obtains that the average rank of $\mathcal{E}$ on $[N,2N]$ is bounded by
\begin{equation}
    \text{Ave\,Rank}_{[N,2N]}(\mathcal{{E}}) \ \leq \ \frac{1}{\sigma} + r + \frac{1}{2} + \left( \frac{0.986}{\sigma} - \frac{2.966}{\sigma^2 \log R} \right) \frac{m_{\mathcal{E}}}{\log R}.
\end{equation}
We use this to estimate the increase in the average rank coming from the $-m_{\mathcal{E}}p$ term as a result of a positive bias. As in \cite{Mil05}, we assume that the one-level density is known for $\sigma = 1$. Fix $m_{\mathcal{E}} = m$ some integer. We get that the lower-order correction is $0.03$ for conductors of size $10^{12}$ with this support. This yields the upper bound
\begin{equation}
    1 + r + \frac{1}{2} + 0.03m \ = \ 1+0.03m + r + \frac{1}{2}.
\end{equation}
Specializing to $m = 1$ gives the upper bound $1.03 + r + 1/2$ on the average rank. This is considerably higher than the observed bound $0.40 + r + 1/2$ due to Fermigier \cite{Fe96}.


\end{document}